\documentclass[12pt,twoside]{amsart}
\usepackage{amssymb}

\numberwithin{equation}{section}
\newtheorem{theorem}{Theorem}[section]

\newtheorem{corollary}[theorem]{Corollary}

\theoremstyle{definition}
\newtheorem{definition}[theorem]{Definition} 

\newtheorem{example}[theorem]{Example}

\newtheorem{notation}[theorem]{Notation}

\begin{document}


\newcommand{\m}[1]{\marginpar{\addtolength{\baselineskip}{-3pt}{\footnotesize
\it #1}}}
\newcommand{\A}{\mathcal{A}}
\newcommand{\K}{\mathcal{K}}
\newcommand{\knd}{\mathcal{K}^{[d]}_n}
\newcommand{\F}{\mathcal{F}}
\newcommand{\N}{\mathbb{N}}
\newcommand{\Z}{\mathbb{Z}}
\newcommand{\pr}{\mathbb{P}}
\newcommand{\I}{\mathcal{I}}
\newcommand{\G}{\mathcal{G}}
\newcommand{\lcm}{\operatorname{lcm}}
\newcommand{\ndp}{N_{d,p}}
\newcommand{\tor}{\operatorname{Tor}}
\newcommand{\reg}{\operatorname{reg}}
\newcommand{\mf}{\mathfrak{m}}
\newcommand{\LL}{\mathcal{L}}
\def\aa{{{\rm \bf a}}}
\def\bb{{{\rm \bf b}}}
\def\cc{{{\rm \bf c}}}
\def\0{{\bf 0}}
\title{Whiskers and sequentially Cohen-Macaulay graphs}
\author{Christopher A. Francisco}
\address{Department of Mathematics, Oklahoma State University, 401 Mathematical Sciences, Stillwater, OK 74078}
\email{chrisf@math.missouri.edu}
\urladdr{http://www.math.missouri.edu/$\sim$chrisf}

\author{Huy T\`ai H\`a}
\address{Department of Mathematics, Tulane University, 6823 St. Charles Ave., New Orleans, LA 70118}
\email{tai@math.tulane.edu}
\urladdr{http://www.math.tulane.edu/$\sim$tai/}
\thanks{The second author is partially supported by the Louisiana Board of Regents Enhancement Grant.}

\subjclass[2000]{}

\begin{abstract}
We investigate how to modify a simple graph $G$ combinatorially to obtain a sequentially Cohen-Macaulay graph. We focus on adding configurations of whiskers to $G$, where to add a whisker one adds a new vertex and an edge connecting this vertex to an existing vertex of $G$. We give various sufficient conditions and necessary conditions on a subset $S$ of the vertices of $G$ so that the graph $G \cup W(S)$, obtained from $G$ by adding a whisker to each vertex in $S$, is a sequentially Cohen-Macaulay graph. For instance, we show that if $S$ is a vertex cover of $G$, then $G \cup W(S)$ is a sequentially Cohen-Macaulay graph. On the other hand, we show that if $G \backslash S$ is not sequentially Cohen-Macaulay, then $G \cup W(S)$ is not a sequentially Cohen-Macaulay graph. Our work is inspired by and generalizes a result of Villarreal on the use of whiskers to get Cohen-Macaulay graphs.
\end{abstract}

\maketitle
\section{Introduction} \label{s.intro}
Let $G = (V_G,E_G)$ be a simple graph (no loops or multiple edges) with vertex set $V_G = \{x_1, \dots, x_n\}$ and edge set $E_G$. We can associate to $G$ a square-free monomial ideal 
\[\I(G) = (x_ix_j ~|~ \{x_i, x_j\} \in E_G) \subset R = k[x_1, \dots, x_n]\] 
The ideal $\I(G)$ is usually referred to as the \textbf{edge ideal} of $G$.

In recent years there has been a flurry of work investigating how the combinatorial data of $G$ appears in algebraic invariants and properties of $\I(G)$. We mention, for example, the works \cite{eghp, FVTchordal, HVT, HVTsurvey, HHZ,  J, JK, K, SVV, V1, Z}. In this paper, we examine how a particular structure of $G$ affects the Cohen-Macaulayness and sequentially Cohen-Macaulayness of its edge ideal. The property of being sequentially Cohen-Macaulay was first introduced by Stanley \cite{Stanley} as a generalization of the well-known property of being Cohen-Macaulay. We recall the definition of sequentially Cohen-Macaulay modules over a polynomial ring.

\begin{definition} \label{d.scm}
Let $M$ be a graded module over $R = k[x_1,\dots,x_n]$. We say that $M$ is \textbf{sequentially Cohen-Macaulay} if there exists a filtration \[ 0=M_0 \subset M_1 \subset \cdots \subset M_r = M \] of $M$ by graded $R$-modules such that $\dim M_i/M_{i-1} < \dim M_{i+1}/M_i$ for all $i$, where $\dim$ denotes Krull dimension, and $M_i/M_{i-1}$ is Cohen-Macaulay for all $i$.
\end{definition}

A graph $G$ is said to be \textbf{(sequentially) Cohen-Macaulay} if $R/\I(G)$ is a (sequentially) Cohen-Macaulay $R$-module. Stanley introduced the notion of sequential Cohen-Macaulayness in connection with work of Bj\"orner and Wachs on nonpure shellability. Pure shellable complexes are Cohen-Macaulay, and Stanley identified sequential Cohen-Macaulayness as the appropriate analogue in the nonpure setting; that is, all nonpure shellable complexes are sequentially Cohen-Macaulay. The notion of sequentially Cohen-Macaulayness has arisen in a number of interesting contexts. For example, Peskine characterized the sequentially Cohen-Macaulay $R$-modules in terms of vanishing or Cohen-Macaulayness of certain Ext modules. For a proof, see Herzog and Sbarra's paper \cite{HSb}, where they also show that if $I$ is a homogeneous ideal, $R/I$ is sequentially Cohen-Macaulay if and only if the local cohomology modules of $R/I$ and $R/\text{Gin}(I)$ have the same Hilbert functions (using the reverse-lex gin). On the combinatorial side, one particularly interesting result is due to Duval, who showed in \cite{D} that algebraic shifting preserves the $h$-triangle of a simplicial complex $\Delta$ if and only if $\Delta$ is sequentially Cohen-Macaulay.

Classifying all Cohen-Macaulay or sequentially Cohen-Macaulay graphs is intractable, and thus it is natural to study some special classes of graphs. Of particular interest is the class of {\it trees} and {\it forests}, or slightly more generally, {\it chordal} graphs. For example, Faridi \cite{Faridi} showed that simplicial trees are sequentially Cohen-Macaulay; the first author and Van Tuyl \cite{FVTchordal} extended this property in the case of graphs to the class of all chordal graphs; and, on the other hand, Herzog, Hibi, and Zheng \cite{HHZ} proved that a chordal graph is Cohen-Macaulay if and only if it is {\it unmixed}.

Our paper complements this work, asking: Given an arbitrary graph $G$, how can one modify $G$ to obtain a sequentially Cohen-Macaulay graph? Motivated by Villarreal's work in \cite{V1}, we investigate the effect of adding ``whiskers'' to a graph. To add a \textbf{whisker} at a vertex $y$ to $G$, one adds a new vertex $x$ and the edge connecting $x$ and $y$ to $G$. We denote by $G \cup W(y)$ the graph obtained from $G$ by adding a whisker at $y$. More generally, if $S \subset V_G$ is a subset of the vertices of $G$, then we denote by $G \cup W(S)$ the graph obtained from $G$ by adding a whisker at each vertex in $S$. The origin of the name is clear from a picture of a whisker added to each vertex of a cycle, and the terminology appears in \cite{SVV}. 

A primary inspiration for this paper is Villarreal's theorem from \cite{V1} (where contributions of Vasconcelos, Herzog, and Fr\"oberg are also cited); see also \cite[Theorem 2.1]{SVV}. He showed that if $G$ is a graph, and $H$ is the graph formed by adding a whisker to {\it every} vertex of $G$, then $H$ is Cohen-Macaulay. This result is sharp: One needs in general to add a whisker to all vertices of $G$, and adding fewer whiskers can actually make a Cohen-Macaulay graph no longer Cohen-Macaulay. 

The goal of this paper is to explore how adding different configurations of whiskers to a graph $G$ affects the weaker property of a graph being sequentially Cohen-Macaulay. 
Our first main result is:

\medskip

\noindent {\bf Theorem~\ref{t.whiskerscm}.} Let $G$ be a simple graph and let $S \subset V_G$. Suppose that $G \backslash S$ is a chordal graph or a five-cycle $C_5$. Then $G \cup W(S)$ is a sequentially Cohen-Macaulay graph.

\medskip

Theorem \ref{t.whiskerscm} has a number of interesting consequences. For example, Corollary \ref{c.manycor} says that if $S \subset V_G$ is a vertex cover, then $G \cup W(S)$ is sequentially Cohen-Macaulay. Thus to create sequentially Cohen-Macaulay graphs by adding whiskers, the number of whiskers is not as important as their configuration. On the other hand, Corollary \ref{c.n-3} says that if $|S| \ge |V_G|-3$, then $G \cup W(S)$ is sequentially Cohen-Macaulay. This gives a bound on the number of vertices so that adding this many whiskers, regardless of their configuration, always results in a sequentially Cohen-Macaulay graph. Furthermore, we recover Villarreal's theorem on creating Cohen-Macaulay graphs in Corollary~\ref{c.whiskercm}.

Our approach uses Alexander duality of edge ideals. If $\Delta$ is a simplicial complex, the faces of the \textbf{Alexander dual} complex $\Delta^*$ are the complements of the nonfaces of $\Delta$. One can define Alexander duality for square-free monomial ideals without reference to simplicial complexes; the duality maps generators of $I$ to primary components of $I^{\vee}$. 

A powerful feature of Alexander duality is that it allows us to link the sequentially Cohen-Macaulay property with homological features of the dual. To do this, Herzog and Hibi introduced the concept of componentwise linearity \cite{HH}. If $I$ is a homogeneous ideal, write $(I_d)$ for the ideal generated by all homogeneous degree $d$ elements of $I$. The ideal $I$ is \textbf{componentwise linear} if for all $d \in \N$, $(I_d)$ has a linear resolution.

All ideals with linear resolutions are componentwise linear, but so are all stable ideals and a number of others; see, e.g., \cite{FVT}. When $I$ is a square-free monomial ideal, Herzog and Hibi found a useful criterion for $I$ to be componentwise linear. Write $(I_{[d]})$ for the ideal generated by all square-free monomials of degree $d$ in $I$. Then $I$ is componentwise linear if and only if for all $d \in \N$, $(I_{[d]})$ has a linear resolution. Moreover, in \cite{HH}, Herzog and Hibi generalized the analogous result on Cohen-Macaulayness from \cite[Theorem 3]{ER}.

\begin{theorem} \label{t.scmcwl}
Let $I$ be a square-free monomial ideal in $R=k[x_1,\dots,x_n]$. $R/I$ is sequentially Cohen-Macaulay over $k$ if and only if $I^{\vee}$ is componentwise linear in $R$.
\end{theorem}

Theorem~\ref{t.scmcwl} allows us to investigate the sequentially Cohen-Macaulay property of an ideal by determining when the Alexander dual is componentwise linear, and Herzog and Takayama's theory of linear quotients \cite{HT} gives a useful way to prove that an ideal is componentwise linear. A monomial ideal $I$ is said to have \textbf{linear quotients} if $I$ has a system of minimal generators $\{u_1, \dots, u_r\}$ with $\deg u_1 \le \dots \le \deg u_r$ such that for all $1 \le i \le r-1$, $(u_1,\dots,u_i):(u_{i+1})$ is generated by linear forms. It is easy to see that if $I$ is an ideal generated in a single degree, and $I$ has linear quotients, then $I$ has a linear resolution. We shall use this observation frequently: If $(I_d)$ has linear quotients for all $d$, then $I$ is componentwise linear. In particular, if $I$ is a square-free monomial ideal, and for all $d$, $(I_{[d]})$ has a linear resolution, then $I$ is componentwise linear. Note that having linear quotients is independent of the characteristic of the field $k$.

The proof of Theorem~\ref{t.whiskerscm} is based upon examining the degree $d$ pieces of the Alexander duals of $\I(G \backslash S)$ and $\I(G \cup W(S))$; we make the following definition:

\begin{definition} \label{d.graphlq}
We say that a graph $G$ has \textbf{dual linear quotients} if for each degree $d \in \N$, $(\I(G)^{\vee}_{[d]})$ has linear quotients.
\end{definition}

If $G$ has dual linear quotients, then $G$ is sequentially Cohen-Macaulay (over a field of any characteristic); this approach has been used in papers of Faridi \cite{Faridi} and the first author and Van Tuyl \cite{FVTchordal}, and we exploit it throughout the paper.

Our second main result addresses the converse of Theorem~\ref{t.whiskerscm}. We give a necessary condition on $G \backslash S$ for $G \cup W(S)$ to be sequentially Cohen-Macaulay.

\medskip

\noindent{\bf Theorem~\ref{t.notscm}.} Let $G$ be a simple graph and let $S \subset V_G$. If $G \backslash S$ is not a sequentially Cohen-Macaulay graph, then $G \cup W(S)$ is not sequentially Cohen-Macaulay.

\medskip

To prove Theorem \ref{t.notscm}, we examine syzygies of the Alexander dual of $\I(G \cup W(S))$ via simplicial homology and {\it upper Koszul simplicial complexes} associated to square-free monomial ideals. Our arguments are inspired by \cite{FVTchordal}.

\medskip

\noindent{\bf Acknowledgement:} The authors would like to thank Adam Van Tuyl for many stimulating discussions on the topic.


\section{Preliminaries} \label{s.prelims}

In this section, we fix some notation for graphs and discuss a result from Alexander duality that we shall use in the rest of the paper.

Throughout, $G = (V_G,E_G)$ denotes a simple graph, which is a graph without any loops or multiple edges. Often, we shall simply write $G$ and not specify its vertex and edge sets. For a vertex $x$ of $G$ we use $N(x)$ to denote the set of neighbors of $x$, which are the vertices connected to $x$ by an edge of $G$. An \textbf{induced subgraph} of $G$ is a subgraph $H$ of $G$ with the property that if $\{z_1, z_2\} \subset V_H$ and $z_1z_2 \in E_G$, then $z_1z_2 \in E_H$.

\begin{notation}
Let $G$ be a simple graph, and let $S=\{y_1,\dots,y_n\}$ be a subset of vertices of $G$. By $G \cup W(S)$ we mean the graph with whiskers $x_iy_i$, for each $1 \le i \le n$, attached to $G$. For simplicity, we shall use $\{x_1y_1, \dots, x_ny_n\}$ to denote $W(S)$ in this case. We use $G \backslash \{y_1,\dots,y_r\}$ to mean the subgraph obtained from $G$ by removing the vertices $y_1,\dots,y_r$ and all edges incident to at least one of these vertices.
\end{notation}

Cycles are important in several of our results. A \textbf{cycle} in a simple graph $G$ is an alternating sequence of distinct vertices and edges $C = v_1e_1v_2e_2 \dots v_{n-1}e_{n-1}v_ne_nv_1$ in which the edge $e_i$ connects the vertices $v_i$ and $v_{i+1}$ ($v_{n+1} \cong v_1$) for all $i$. In this case, we say $C$ has \textbf{length} $n$ and call $C$ an $n$-cycle. We shall also use $C_n$ to denote an $n$-cycle.

We are particularly interested in chordal graphs, which includes all trees and forests and has been the object of much study in recent years. We call a graph $G$ \textbf{chordal} if for all $n \ge 4$, every $n$-cycle in $G$ has a \textbf{chord}, which is an edge connecting two non-consecutive vertices of the cycle. Notice that an induced subgraph of a chordal graph is also chordal. Thus, as a byproduct of the arguments of \cite[Theorem 3.2]{FVTchordal}, we get a theorem on chordal graphs that we use in inductive arguments in the next section:

\begin{theorem} \label{t.chordalscm}
Let $G$ be a chordal graph, and let $H$ be an arbitrary induced subgraph of $G$. Then $H$ has dual linear quotients.
\end{theorem}

The generators of the Alexander dual of an edge ideal represent covers of the associated graph. We recall some terminology relating to these covers.

\begin{definition} \label{d.vc}
A \textbf{vertex cover} of a graph $G$ is a subset $V \subset V_G$ of the vertices of $G$ such that every edge of $G$ is incident to at least one vertex of $V$ (in particular, isolated vertices need not appear in a vertex cover). If $V$ is a vertex cover of $G$, then it is a \textbf{minimal vertex cover} of $G$ if no proper subset of $V$ is a vertex cover of $G$. For simplicity, we write vertex covers as monomials, so $\{z_1,\dots,z_r\}$ will be written as $z_1 \cdots z_r$. A graph is said to be \textbf{unmixed} if every minimal vertex cover has the same cardinality.
\end{definition}

Using Alexander duality, one can describe how being Cohen-Macaulay differs from being sequentially Cohen-Macaulay in the square-free monomial ideal case. By \cite[Lemma 3.6]{FVTchordal}, if $I \subset R$ is square-free monomial ideal, then $R/I$ is Cohen-Macaulay if and only if $R/I$ is sequentially Cohen-Macaulay and $I$ is unmixed. Consequently, because the unmixedness of a graph is a combinatorial property (on the cardinality of minimal vertex covers), investigating the Cohen-Macaulayness of a graph reduces to studying the sequentially Cohen-Macaulayness of such a graph. In particular, a sequentially Cohen-Macaulay graph is Cohen-Macaulay if and only if it is unmixed.


\section{Whiskers and sequentially Cohen-Macaulay graphs}\label{s.scm}

In this section, we explore how to add a configuration of whiskers to an arbitrary graph to create a sequentially Cohen-Macaulay graph. Let $G$ be a graph, and let $S \subset V_G$. Our primary question is: What conditions on $S$ make $G \cup W(S)$ sequentially Cohen-Macaulay? Because being sequentially Cohen-Macaulay is a weaker property than being Cohen-Macaulay, one expects that $S$ needs not be all of $V_G$ to ensure that $G \cup W(S)$ is sequentially Cohen-Macaulay. The focus of this section is on sufficient conditions that guarantee the sequential Cohen-Macaulayness of $G \cup W(S)$, and our results also recover Villarreal's theorem on Cohen-Macaulayness as a consequence.

Because the vertex covers of $G$ are the generators of the Alexander dual of $\I(G)$, we are often interested in ways to partition the set of vertex covers of $G$ of a particular cardinality. For any graph $G$ and vertex $x \in V_G$ with $N(x)=\{y_1,\dots,y_t\}$, we can decompose the set of vertex covers of $G$ of size $d$ in the following way: Any vertex cover of $G$ of size $d$ is either $x$ times a vertex cover of $G \backslash \{x\}$ of size $d-1$, or it is $y_1\cdots y_t$ times a vertex cover of $G \backslash \{x,y_1,\dots,y_t\}$ of size $d-t$.
For our purposes, we frequently consider the case in which $G$ contains a whisker $xy$, where $x$ is the vertex of degree one. In this case, the vertex covers of $G$ are decomposed based on covers of $G \backslash \{x\}$ and covers of $G \backslash \{x,y\}$. In particular, the set of vertex covers of $G$ of size $d$ is the union of $x$ times the vertex covers of $G \backslash \{x\}$ of size $d-1$ and $y$ times the vertex covers of $G \backslash \{x,y\}$ of size $d-1$.

The next theorem is the first step in exploring how to add whiskers to a graph to make it sequentially Cohen-Macaulay. 

\begin{theorem} \label{t.subgraphlq}
Let $G'$ be a simple graph and let $S = \{y_1, \dots, y_n\} \subset V_{G'}$ be a subset of the vertices of $G'$. Let $\{x_1y_1, \dots, x_ny_n\}$ be whiskers of $G'$ at $S$ and let $G = G' \cup W(S)$. Suppose that if $H \subset G$ is an induced subgraph of $G$ such that both
\begin{enumerate}
\item[$(i)$] $\{x_1,\dots,x_{n-1}\} \subset V_H$, and
\item[$(ii)$] $x_n \not \in H$ and $y_n \not \in H$
\end{enumerate}
hold, then $H$ has dual linear quotients. Then all induced subgraphs $K \subset G$ such that $\{x_1,\dots,x_n\} \subset V_K$ have dual linear quotients.
\end{theorem}

\begin{proof} For simplicity of notation, let $\{z_1, \dots, z_r\} = V_{G'} \backslash S$. Fix an induced subgraph $K \subset G$ as in the statement of the theorem. Consider first the case in which $y_n \not\in K$. Then $x_n$ is an isolated vertex of $K$. Let $H$ be the graph obtained from $K$ by removing the isolated vertex $x_n$. Clearly, $H$ satisfies properties $(i)$ and $(ii)$ of the hypothesis. Thus, $H$ has dual linear quotients, i.e., $(\I(H)^{\vee}_{[d]})$ has linear quotients for all $d \in \N$. Since the only edge to which $x_n$ is incident in $G$ is $x_ny_n$, the minimal generating set of $\I(K)$ is the same as the minimal generating set of $\I(H)$; thus the minimal generating sets of $\I(K)^{\vee}$ and $\I(H)^{\vee}$ are the same, though the first is an ideal of $k[x_1,\dots,x_n,y_1,\dots,y_{n-1},z_1,\dots,z_r]$, and the second is an ideal of $k[x_1,\dots,x_{n-1},y_1,\dots,y_{n-1},z_1,\dots,z_r]$. Thus for all $d \in \N$, by \cite[Lemma 2.9]{FVT}, since $(\I(H)^{\vee}_{[d]})$ has linear quotients, so does $(\I(K)^{\vee}_{[d]})$.

Consider instead the case in which $y_n \in K$. Let $H$ be the subgraph of $K$ obtained by removing $x_n$, $y_n$, and all edges incident to $x_n$ or $y_n$. Again, $H$ satisfies $(i)$ and $(ii)$ of the hypotheses; and thus, $H$ has dual linear quotients. Fix a degree $d \in \N$. Let $A_1,\dots,A_a$ be the monomials that represent all vertex covers of $K \backslash \{x_n\}$ of size $d-1$, and let $B_1,\dots,B_b$ be the monomials that represent all vertex covers of $H=K \backslash \{x_n,y_n\}$ of size $d-1$; that is, $(B_1,\dots,B_b)=(\I(H)^{\vee}_{[d-1]})$. We have \[ (\I(K)^{\vee}_{[d]}) = x_n(A_1,\dots,A_a) + y_n(B_1,\dots,B_b).\] The $A_i$s are monomials in the variables $x_1,\dots,x_{n-1},y_1,\dots,y_n,z_1,\dots,z_r$, and the $B_i$s are monomials in the variables $x_1,\dots,x_{n-1},y_1,\dots,y_{n-1},z_1,\dots,z_r$. Since $H$ has dual linear quotients, $(\I(H)^{\vee}_{[d]})$ has linear quotients for all $d$. We may assume that the $B_i$s are indexed in the order that gives linear quotients (that is, $(B_1,\dots,B_{i-1}):B_i$ is generated by a subset of the variables for all $i$).

We wish to show that the ideal $(y_nB_1, \dots, y_nB_b, x_nA_1, \dots,x_nA_a)$ has linear quotients. Since $(B_1,\dots,B_b)$ has linear quotients (in that order), it suffices to show that for all $j$, \[ (y_nB_1,\dots, y_nB_b, x_nA_1, \dots, x_nA_{j-1}):x_nA_j \] is generated by a subset of the variables. To this end, we consider two possibilities.

Suppose first that $y_n$ divides $A_j$. Then $x_nA_j=x_ny_nC$, where $C$ is a vertex cover of $H = K \backslash \{x_n,y_n\}$ of size $d-2$. Let $T$ be the set of variables in $\{x_1,\dots,x_{n-1}$, $y_1,\dots,y_{n-1}$, $z_1,\dots,z_r\}$ which are not in the support of $C$, and suppose $u \in T$. Then $uC$ is a vertex cover of $H$ of size $d-1$, so it is one of the $B_i$. Therefore for any $u \in T$, $ux_nA_j \in (y_nB_1,\dots,y_nB_b)$. Moreover, note that $(y_nB_1,\dots,y_nB_b,x_nA_1,\dots,x_nA_{j-1})$ is a square-free monomial ideal; thus, if $m$ is a minimal monomial generator of 
$$(y_nB_1, \dots, y_nB_b, x_nA_1, \dots, x_nA_{j-1}) : x_nA_j,$$ then $m$ is square-free. Hence $x_n$, $y_n$, and variables dividing $C$ do not divide $m$. Thus \[ (y_nB_1,\dots,y_nB_b,x_nA_1,\dots,x_nA_{j-1}):x_nA_j = (\mbox{all variables } u \in T).\]

Next we assume that $y_n$ does not divide $A_j$. Note that any $A_j$ that is not divisible by $y_n$ is one of the $B_i$s because a cover of $K \backslash \{x_n\}$ not containing $y_n$ is a cover of $H$. Thus $A_j=B_{i_j}$ for some $i_j$. Consider a monomial $m$ for which $mx_nA_j \in (y_nB_1,\dots,y_nB_b)$. Then since $y_n$ does not divide $A_j$, $y_n$ must divide $m$. But $y_nx_nA_j=y_nx_nB_{i_j} \in (y_nB_1,\dots,y_nB_b)$, so $y_nx_nA_j \in (y_nB_1,\dots,y_nB_b)$. Hence $(y_nB_1,\dots,y_nB_b):x_nA_j = (y_n)$.

The last remaining situation is when $y_n$ does not divide $A_j$, and $m$ is a monomial such that $mx_nA_j$ lands in the ideal $(x_nA_1, \dots, x_nA_{j-1})$. This case requires a bit more work. We need to specify an order for the $A_i$ monomials. Note that, so far we have not used any feature of the ordering of the $A_i$s. We may pick an order so that all the $A_i$s not divisible by $y_n$ are indexed first, and those that are divisible by $y_n$ are last. Suppose that $\{A_1, \dots, A_t\}$ are all the $A_i$s that are not divisible by $y_n$, and $A_l = B_{i_l}$ for $l = 1, \dots, t$. For each $1 \le l \le t$, $A_l = B_{i_l}$ is a vertex cover of $K \backslash \{x_n\}$ not containing $y_n$, so it is divisible by all variables in $N(y_n) \backslash \{x_n\}$. Let $D$ be the monomial given by the variables in $N(y_n) \backslash \{x_n\}$. For $1 \le l \le t$, let $C_l = A_l/D$. Then, clearly, $\{C_1,\dots, C_t\}$ are the vertex covers of $L = K \backslash \{x_n,y_n,N(y_n)\}$ of size $d-1-u$, where $u = |N(y_n)|-1$. Conversely, if $C$ is a vertex cover of $L$, then $CD$ is a vertex cover of $K \backslash \{x_n\}$ not containing $y_n$. Thus $\{C_1, \dots, C_t\}$ are all vertex covers of $L$ of size $d-1-u$. Since $x_n \in N(y_n)$ but $x_j \not \in N(y_n)$ for $j \not = n$, it is easy to see that $L$ satisfies $(i)$ and $(ii)$ of the hypotheses. Therefore, $L$ has dual linear quotients. This implies that the ideal $(C_1, \dots, C_t)$ has linear quotients. We shall reindex $\{A_1, \dots, A_t\}$ so that $C_1, \dots, C_t$ is the order of the generators in which $(C_1, \dots, C_t)$ has linear quotients.

Now suppose that $m$ is a monomial so that $mx_nA_j \in (x_nA_1, \dots, x_nA_{j-1})$. Dividing by the monomial given by $N(y_n)$ (including $x_n$), we have $m C_j \in (C_1,\dots,C_{j-1})$. Since $(C_1,\dots,C_t)$ has linear quotients, $(C_1,\dots,C_{j-1}):C_j=(x_{p_1},\dots,x_{p_v})$ for some subset of the variables. Thus if $mx_nA_j \in (x_nA_1, \dots, x_nA_{j-1})$, then some variable $x_{p_w}$ divides $m$.

We have shown that $(y_nB_1,\dots,y_nB_b, x_nA_1,\dots, x_nA_a)$ has linear quotients. This is true for any $d \in \N$. Hence, the conclusion follows.
\end{proof}

We are now ready to prove our first main result.

\begin{theorem} \label{t.whiskerchordal}
Let $G$ be a simple graph. Let $S \subset V_G$ be such that $G \backslash S$ is a chordal graph. Then $G \cup W(S)$ is a sequentially Cohen-Macaulay graph.
\end{theorem}

\begin{proof} Let $S = \{y_1, \dots, y_n\}$ and $W(S) = \{x_1y_1, \dots, x_ny_n\}$. It suffices to show that $G \cup W(S)$ has dual linear quotients.

We shall first construct a class of subgraphs of $G$ as follows. Let $G_0=G \backslash S$. Let $G_1 = G_0 \cup \{x_1,y_1\} \cup \{ \mbox{edges of $G$ incident to vertices of $V_{G_0} \cup \{x_1, y_1\}$}\}$. More generally, for $1 \le i \le n$, let $G_i = G_{i-1} \cup \{x_i, y_i\} \cup \{ \mbox{edges of $G$ incident to vertices of $V_{G_{i-1}} \cup \{x_i, y_i\}$}\}$. Observe that $G_n = G$. Now, the conclusion will follow if we can show that every induced subgraph $K$ of $G$ (in particular, $G$ itself) containing $\{x_1, \dots, x_n\}$ has dual linear quotients. To this end, we shall use induction on $i$ to show that every induced subgraph $K$ of $G_i$ containing $\{x_1, \dots, x_i\}$ has dual linear quotients for $i = 0, \dots, n$.

Indeed, for $i = 0$, the assertion follows from Theorem \ref{t.chordalscm}. Suppose $i \ge 1$. Consider an arbitrary induced subgraph $H$ of $G_i$ such that $\{x_1, \dots, x_{i-1}\} \subset V_H$, $x_i \not\in V_H$ and $y_i \not\in V_H$. Then $H$ is also an induced subgraph of $G_{i-1} = G_i \backslash \{x_i, y_i\}$. Thus, by induction, $H$ has dual linear quotients. It now follows from Theorem \ref{t.subgraphlq} that every induced subgraph $K$ of $G_i$ with $\{x_1, \dots, x_i\} \subset V_K$ has dual linear quotients. The theorem is proved.
\end{proof}




We can extend Theorem \ref{t.whiskerchordal} slightly.

\begin{theorem} \label{t.whiskerscm}
Let $G$ be a simple graph and let $S \subset V_G$. Suppose $G \backslash S$ is a chordal graph or a five-cycle $C_5$. Then $G \cup W(S)$ is a sequentially Cohen-Macaulay graph.
\end{theorem}

\begin{proof} If $G \backslash S$ is chordal, then the conclusion is Theorem \ref{t.whiskerchordal}. Suppose that $G \backslash S = C_5$. The inductive argument of Theorem \ref{t.whiskerchordal} rests on the fact that every induced subgraph of $G_0 = G \backslash S$ has dual linear quotients. In our current situation, $G_0 = G \backslash S = C_5$. Moreover, $\I(C_5)^{\vee}=(x_1x_2x_4, x_1x_3x_4,x_1x_3x_5,x_2x_3x_5,x_2x_4x_5)$ has dual linear quotients in the given order of the generators. Any proper subgraph of $C_5$ is a forest and thus also has dual linear quotients. Hence, arguments similar to those in Theorem \ref{t.whiskerchordal} yield the assertion.
\end{proof}

Theorems~\ref{t.whiskerchordal} and \ref{t.whiskerscm} give many interesting corollaries about configurations of whiskers that can be added to obtain a sequentially Cohen-Macaulay graph. Part (i) below gives a particularly easy way to ensure one constructs a sequentially Cohen-Macaulay graph.

\begin{corollary} \label{c.manycor}
Let $G$ be a simple graph, and let $S \subset V_G$. 
\begin{enumerate}
\item[(i)] If $S$ is a vertex cover of $G$, then $G \cup W(S)$ is sequentially Cohen-Macaulay.
\item[(ii)] If $G \backslash S$ is a forest, then $G \cup W(S)$ is sequentially Cohen-Macaulay.
\item[(iii)] If $G=C$ is a cycle, and $y$ is a vertex of $C$, then $C \cup W(y)$ is sequentially Cohen-Macaulay.
\end{enumerate}
\end{corollary}

\begin{proof} For (i), since $S$ is a vertex cover of $G$, $G \backslash S$ is a graph of isolated vertices. A graph without any edges is clearly a chordal graph. Thus the assertion is a consequence of Theorem \ref{t.whiskerchordal}. (ii) is immediate from Theorem~\ref{t.whiskerchordal} since every forest is a chordal graph. Finally, removing a vertex from a cycle produces a tree, so (iii) follows from (ii).
\end{proof}

Corollary~\ref{c.manycor}(iii) allows one to make a cycle into a sequentially Cohen-Macaulay graph quite easily; only one whisker is necessary. With Van Tuyl, we noticed this phenomenon after doing many computations in the computer algebra system Macaulay 2 \cite{M2}, and it was a primary initial motivation for this paper.

Notice that Corollary~\ref{c.manycor} states that to obtain sequentially Cohen-Macaulay graphs, the number of whiskers is not as important as their configuration. Our next corollary complements Corollary \ref{c.manycor} to give a bound on the number of whiskers to add to a graph, regardless of how they are picked, to obtain a sequentially Cohen-Macaulay graph.

\begin{corollary} \label{c.n-3}
Let $G$ be a simple graph and let $S \subset V_G$. Assume that $|S| \ge |V_G| - 3$. Then $G \cup W(S)$ is a sequentially Cohen-Macaulay graph.
\end{corollary}

\begin{proof} Since $|S| \ge |V_G| - 3$, $G \backslash S$ is a graph on at most 3 vertices. Thus, $G \backslash S$ is either a three-cycle, a tree, or a graph with isolated vertices. These are all chordal graphs, and hence, the conclusion follows from Theorem \ref{t.whiskerchordal}.
\end{proof}

We shall see in Example \ref{e.allbutfour} that the bound $|V_G|-3$ in Corollary \ref{c.n-3} is sharp. Additionally, Corollary \ref{c.n-3} allows us to recover Villarreal's theorem.

\begin{corollary} \label{c.whiskercm}
Let $G$ be a simple graph with vertex set $V_G$. Then $G \cup W(V_G)$ is a Cohen-Macaulay graph.
\end{corollary}

\begin{proof} By Corollary \ref{c.n-3} we know that $G' = G \cup W(V_G)$ is sequentially Cohen-Macaulay. Thus it suffices to show that $G'$ is unmixed; i.e., all minimal vertex covers of $G'$ have the same cardinality. Suppose $V_G = \{y_1, \dots, y_n\}$ and $W(V_G) = \{x_1y_1, \dots, x_ny_n\}$. Let $V$ be an arbitrary minimal vertex cover of $G'$. Clearly, for each $i=1, \dots, n$, $V$ has to contain one of the vertices $\{x_i,y_i\}$ (to cover the edge $x_iy_i$). Moreover, since $V$ is minimal, for each $i=1, \dots, n$, $V$ contains exactly one of the vertices $\{x_i,y_i\}$. Hence, $|V| = n$. This is true for any minimal vertex cover $V$ of $G'$. Thus, the assertion follows.
\end{proof}

In the final theorem of this section, we isolate the condition from the proof of Theorem~\ref{t.whiskerchordal} that yields that result and its corollaries.

\begin{theorem} \label{t.iff}
Let $G$ be a simple graph with $S$ a subset of its vertex set. Then all induced subgraphs of $G \backslash S$ have dual linear quotients if and only if all induced subgraphs of $G \cup W(S)$ have dual linear quotients.
\end{theorem}

\begin{proof}
If all induced subgraphs of $G \cup W(S)$ have dual linear quotients, then so do all induced subgraphs of $G \backslash S$ since any induced subgraph of $G \backslash S$ is an induced subgraph of $G \cup W(S)$. Assume instead that all induced subgraphs of $G \backslash S$ have dual linear quotients. Then an argument identical to the proof of Theorem~\ref{t.whiskerchordal} shows that all induced subgraphs of $G \cup W(S)$ have dual linear quotients.
\end{proof}

We give two examples showing that the hypotheses of Theorem~\ref{t.iff} cannot easily be weakened.

\begin{example} \label{e.notdlq}
Let $G$ be a simple graph with $S \subset V_G$, and assume that $G \backslash S$ has dual linear quotients. In this example, we show that if there exists a subgraph of $G \backslash S$ that does not have dual linear quotients, then $G \cup W(S)$ may fail to be sequentially Cohen-Macaulay.

Let $G$ be the graph on the vertex set $V_G=\{x_1,\dots,x_6\}$ together with edge set $E_G=\{x_1x_2,x_2x_3,x_3x_4,x_1x_4,x_3x_5,x_4x_5,x_5x_6\}$. Let $S=\{x_6\}$, so $G \cup W(S)$ is the graph $G$ along with a new vertex $x_7$ and edge $x_6x_7$. Then \[ \I(G\backslash S)^{\vee} = (x_1x_3x_4,x_2x_3x_4,x_1x_3x_5,x_2x_4x_5),\] which has linear quotients in the order in which the generators are listed. Hence the graph $G \backslash S$ has dual linear quotients. Note, however, that not all induced subgraphs of $G \backslash S$ have dual linear quotients; the four-cycle comprised of the vertices $\{x_1,\dots,x_4\}$ is not even sequentially Cohen-Macaulay.

Now we consider $G \cup W(S)$. We have \[ \I(G \cup W(S))^{\vee} = (x_1x_3x_4x_6, x_2x_3x_4x_6, x_1x_3x_5x_6, x_2x_4x_5x_6, x_1x_3x_5x_7, x_2x_4x_5x_7). \] The minimal graded free resolution of $\I(G \cup W(S))^{\vee}$ is \[ 0 \longrightarrow R(-7) \longrightarrow R(-5)^5 \oplus R(-6) \longrightarrow R(-4)^6 \longrightarrow \I(G \cup W(S))^{\vee} \longrightarrow 0,\] where $R=k[x_1,\dots,x_7]$. Therefore $\I(G \cup W(S))^{\vee}$ is not componentwise linear because of the syzygies in degrees six and seven. Hence $G \cup W(S)$ is not sequentially Cohen-Macaulay.
\end{example}

\begin{example} \label{e.dlq}
Again we assume that $G$ is a simple graph with $S \subset V_G$ such that $G \backslash S$ has dual linear quotients. Now we show that even if there exists a subgraph of $G \backslash S$ that does not have dual linear quotients (and, in fact, is not sequentially Cohen-Macaulay), $G \cup W(S)$ may itself have dual linear quotients.

Let $G$ be the graph with $V_G=\{x_1,\dots,x_6\}$ and edge set \[E_G=\{x_1x_2,x_2x_3,x_3x_4,x_1x_4,x_3x_5,x_4x_5,x_2x_6,x_3x_6\}.\] Let $S=\{x_6\}$. Then $\I(G \backslash S)^{\vee}$ is the same as in Example~\ref{e.notdlq}, and hence $G \backslash S$ has dual linear quotients. Note that the induced subgraph on the vertices $\{x_1,\dots,x_4\}$ is a four-cycle, which is not sequentially Cohen-Macaulay.

We claim that $G \cup W(S)$ has dual linear quotients. The dual of the edge ideal is \[\I(G \cup W(S))^{\vee} = (x_1x_3x_4x_6, x_2x_3x_4x_6, x_1x_3x_5x_6, x_2x_4x_5x_6, x_2x_3x_4x_7, x_1x_2x_3x_5x_7).\] One can check that this ideal has linear quotients with respect to the order in which the generators are listed. Therefore $G \cup W(S)$ is sequentially Cohen-Macaulay.
\end{example}

Consequently, if $G \backslash S$ has dual linear quotients, but there exists a subgraph of $G \backslash S$ without dual linear quotients, then $G \cup W(S)$ may or may not be sequentially Cohen-Macaulay. This is the primary reason for our techniques in the proofs in Section~\ref{s.scm}; in our inductive approach, we assume that {\it all} induced subgraphs of a certain type have dual linear quotients to avoid cases like Example~\ref{e.notdlq}.

We conclude this section by remarking that it is difficult to find results analogous to Theorem~\ref{t.iff} for Cohen-Macaulay graphs. One is tempted to conjecture that if $G$ is a simple graph, and $S$ is a subset of $V_G$ such that all induced subgraphs of $G\backslash S$ have dual linear quotients and are unmixed, then $G \cup W(S)$ is Cohen-Macaulay. Unfortunately, this is false. An easy counterexample is the case in which $G$ is one edge connecting two vertices $y_1$ and $y_2$. $G \backslash \{y_1\}$ trivially has dual linear quotients and is unmixed. However, $k[x,y_1,y_2]/(xy_1,y_1y_2)$ is not Cohen-Macaulay. The difficulty in searching for the appropriate analogue to Theorem~\ref{t.iff} is guaranteeing the unmixedness of $G \cup W(S)$.


\section{Whiskers and non-sequentially Cohen-Macaulay graphs} \label{s.notscm}

In the previous section, we gave sufficient conditions for getting sequentially Cohen-Macaulay graphs by adding whiskers. In this section, we address the converse problem. Our primary interest is necessary conditions on a graph $G$ and $S \subset V_G$ so that $G \cup W(S)$ has a chance to be sequentially Cohen-Macaulay.

To show that certain graphs are not sequentially Cohen-Macaulay, we exploit Alexander duality and show that the dual of the edge ideal is not componentwise linear. This requires investigating the syzygies of the dual, and to do that, we use simplicial homology.

Define a square-free vector to be a vector with its entries in $\{0,1\}$. For any monomial ideal $M$, we define the \textbf{upper Koszul simplicial complex of} {\boldmath $M$}:
\[K^\bb(M)=\, \, \mbox{\{square-free vectors \aa \, such that } \frac{x^\bb}{x^\aa} \in M\}.\]
See, e.g., \cite{MS}. Using the relation \[ \beta_{i,\bb}(M) = \dim_k \tilde{H}_{i-1}(K^\bb(M),k), \] which is \cite[Theorem 1.34]{MS}, we can compute the $\mathbb N^n$-graded Betti numbers of $M$. We use this technique in the following theorem.

\begin{theorem} \label{t.notscm}
Let $G$ be a simple graph and let $S \subset V_G$. If $G \backslash S$ is not sequentially Cohen-Macaulay, then $G \cup W(S)$ is not sequentially Cohen-Macaulay.
\end{theorem}

\begin{proof} Let $\{z_1, \dots, z_r\} = V_G \backslash S$. Because $G \backslash S$ is not sequentially Cohen-Macaulay, there exists $d$ such that $I:=(\I(G \backslash S)^{\vee}_{[d]})$ has a nonlinear $i$th syzygy in its minimal free resolution. Suppose this nonlinear syzygy occurs in the multi-degree $\bb$ corresponding to the square-free monomial $z_{i_1} \dots z_{i_l}$ for some $l > d+i$. Then $\dim_k \tilde{H}_{i-1}(K^{\bb}(I),k) \not = 0$.

Let $\{y_1, \dots, y_n\} = S$, and let $W(S) = \{x_1y_1, \dots, x_ny_n\}$. Let $J=(\I(G \cup W(S))^{\vee}_{[d+n]})$. Let $\cc$ be the multi-degree of the monomial $z_{i_1} \dots z_{i_l} y_1 \dots y_n$. We claim that the simplicial complexes $K^{\bb}(I)$ and $K^{\cc}(J)$ are the same. By definition, a square-free vector $\aa$ is in $K^{\cc}(J)$ if and only if $\frac{x^\cc}{x^\aa} \in (\I(G \cup W(S))^\vee_{[d+n]})$. In other words, a square-free vector $\aa$ is in $K^\cc(J)$ if and only if the square-free monomial corresponding to $\cc-\aa$ gives a vertex cover of $G \cup W(S)$. Since $\cc$ has 0 in entries corresponding to $\{x_1, \dots, x_n\}$, if $\cc-\aa$ gives a vertex cover of $G \cup W(S)$, then $\cc-\aa$ must have 1 in all entries corresponding to $\{y_1, \dots, y_n\}$, and $\aa$ must have 0 in all entries corresponding to $\{x_1, \dots, x_n\}$. Therefore, the only places in which $\aa$ may be nonzero are in entries corresponding to the $z_i$s. Thus the vectors $\aa$ such that $\cc-\aa$ gives a vertex cover of $G \cup W(S)$ are exactly the same as the vectors $(\aa', \0)$, where $\0 = (0, \dots, 0) \in \Z^{2n}$ appears in the entries corresponding to $\{x_1, \dots, x_n,y_1,\dots,y_n\}$, so that $\bb-\aa'$ gives a vertex cover of $G \backslash S$. Hence the $\aa'$ in $K^{\cc}(J)$ are exactly the $\aa$ in $K^{\bb}(I)$. This implies that the simplicial complexes $K^\bb(I)$ and $K^\cc(J)$ are the same.

We now have $\dim_k \tilde{H}_{i-1}(K^{\cc}(J),k) = \dim_k \tilde{H}_{i-1}(K^\bb(I),k) \not = 0$. This gives a nonlinear $i$th syzygy of $J$ since $l> d+i$ implies $l+n>i+d+n$, and $J$ is generated in degree $d+n$. Hence, $J$ does not have a linear resolution, and thus $G \cup W(S)$ is not sequentially Cohen-Macaulay by Theorem \ref{t.scmcwl}.
\end{proof}

As a corollary, we can identify certain vertex sets to which adding whiskers does not yield a sequentially Cohen-Macaulay graph.

\begin{corollary} \label{c.cyclenotscm}
Let $G$ be a simple graph and let $S \subset V_G$. Suppose that $G \backslash S = C_n$, where $n$ is neither three nor five. Then $G \cup W(S)$ is not sequentially Cohen-Macaulay.
\end{corollary}

\begin{proof} Since $n$ is neither three nor five, by \cite[Proposition 4.1]{FVTchordal}, $C_n$ is not sequentially Cohen-Macaulay. The assertion is a consequence of Theorem \ref{t.notscm}.
\end{proof}

Corollary~\ref{c.cyclenotscm} helps show the sharpness of Corollary \ref{c.n-3}, as illustrated in the following example. 

\begin{example} \label{e.allbutfour}
Let $G$ be a four-cycle with vertices $\{y_1, \dots, y_4\}$ and an edge $y_1y$ attached at $y_1$, so $\I(G)=(y_1y_2,y_2y_3,y_3y_4,y_1y_4,y_1y).$ Let $S=\{y\}$, and add a whisker $xy$ to $G$ to form $G \cup W(S)$. Then \[\I(G \cup W(S))^{\vee} = (y_1y_3y,y_2y_4y,y_1y_3x,y_1y_2y_4x) \subset R=k[y_1,\dots,y_4, y,x].\] Note that $(\I(G \cup W(S))_3^{\vee})$ has minimal graded free resolution \[ 0 \longrightarrow R(-4) \oplus R(-5) \longrightarrow R(-3)^3 \longrightarrow I \longrightarrow 0,\] which is not a linear resolution. Hence the Alexander dual of $\I(G \cup W(S))$ is not componentwise linear, and therefore $G \cup W(S)$ is not sequentially Cohen-Macaulay.
\end{example}

The primary case we have not considered is when a graph is sequentially Cohen-Macaulay but does not have dual linear quotients. This case will require different methods; sequential Cohen-Macaulayness can depend on the underlying field $k$, but having dual linear quotients is independent of the field. We give an example of this phenomenon.

\begin{example} \label{e.chardepend}
Let $\Delta$ be a minimal triangulation of the real projective plane. Then $\Delta$ is Cohen-Macaulay over a field $k$ if and only if the characteristic of $k$ is not two \cite{BH}. From this, we can construct an example of a graph $G$ that is Cohen-Macaulay if and only if the base field does not have characteristic two, using the method described in \cite{HHZ}. Let $P$ be the face poset of $\Delta$; then the order complex of $P$ has the property that all its minimal nonfaces are subsets of cardinality two, so the associated Stanley-Reisner ideal is generated by degree two monomials and hence is an edge ideal (with a large number of generators). The polynomial ring modulo this edge ideal is Cohen-Macaulay if and only if the ground field does not have characteristic two, just like the simplicial complex $\Delta$.
\end{example}

We know of no example of a graph $G$ with a small number of vertices that is (sequentially) Cohen-Macaulay over fields of some characteristics but not others, and it would be interesting to know of small examples if they exist.



\begin{thebibliography}{99}
\bibitem{BH} W.~Bruns and J.~Herzog, \emph{Cohen-{M}acaulay rings}, in \emph{Cambridge Studies in Advanced Mathematics} \textbf{39}, Cambridge University Press, Cambridge, 1993.

\bibitem{D} A. Duval, Algebraic shifting and sequentially Cohen-Macaulay simplicial complexes. \emph{Electron. J. Combin.} {\bf 3} (1996), no. 1, Research Paper 21.

\bibitem{eghp} D. Eisenbud, M. Green, K. Hulek and S. Popescu,
Restricting linear syzygies: algebra and geometry. \emph{Compos. Math.} {\bf 141} (2005), no. 6, 1460--1478.

\bibitem{ER} J. Eagon and V. Reiner, Resolutions of Stanley-Reisner rings and Alexander duality. \emph{J. Pure Appl. Algebra} {\bf 130} (1998), no. 3, 265--275.

\bibitem{Faridi} S. Faridi, Simplicial trees are sequentially Cohen-Macaulay. \emph{J. Pure Appl. Algebra} {\bf 190} (2003), 121--136.

\bibitem{FVTchordal} C. Francisco and A. Van Tuyl, Sequentially Cohen-Macaulay edge ideals. \emph{Proc. Amer. Math. Soc.} {\bf 135} (2007), no. 8, 2327--2337.

\bibitem{FVT} C. Francisco and A. Van Tuyl, Some families of componentwise linear monomial ideals. To appear, \emph{Nagoya Math. J.} {\tt math.AC/0508589}

\bibitem{HVT} Huy T\`ai H\`a, A. Van Tuyl, Splittable ideals and the resolutions
of monomial ideals. \emph{J. Algebra} {\bf 309} (2007), 405--425.

\bibitem{HVTsurvey} Huy T\`ai H\`a, A. Van Tuyl, Resolutions of square-free monomial ideals via facet ideals: a survey. To appear, \emph{Contemp. Math.} {\tt math.AC/0604301}

\bibitem{HH} J. Herzog and T. Hibi, Componentwise linear ideals.  \emph{Nagoya Math. J.}  {\bf 153} (1999), 141--153.

\bibitem{HHZ} J. Herzog, T. Hibi, and X. Zheng, Cohen-Macaulay chordal graphs. \emph{J. Combin. Theory Ser. A} {\bf 113} (2006), no. 5, 911--916.

\bibitem{HSb} J. Herzog and E. Sbarra, Sequentially Cohen-Macaulay modules and local cohomology. \emph{Algebra, arithmetic and geometry, Part I, II (Mumbai, 2000), Tata Inst. Fund. Res. Stud. Math.} {\bf 16} (2002), 327--340.

\bibitem{HT} J. Herzog and Y. Takayama, Resolutions by mapping cones. The Roos Festschrift volume, 2.  \emph{Homology Homotopy Appl.}  {\bf 4}  (2002),  no. 2, part 2, 277--294 (electronic).

\bibitem{J} S. Jacques, Betti numbers of graph ideals. Ph.D. Thesis,
University of Sheffield, 2004. {\tt math.AC/0410107}

\bibitem{JK} S. Jacques, M. Katzman, The Betti numbers of forests.
(2005) Preprint. {\tt math.AC/0401226}

\bibitem{K} M. Katzman, Characteristic-independence of Betti numbers
of graph ideals. \emph{J. Combin. Theory, Ser. A} {\bf 113} (2006), 435--454.

\bibitem{M2}
D.~R. Grayson and M.~E. Stillman, \emph{Macaulay 2, a software system for
  research in algebraic geometry}.
\newblock \verb|http://www.math.uiuc.edu/Macaulay2/|.

\bibitem{MS} E. Miller and B. Sturmfels, {\it Combinatorial commutative algebra}. Springer, 2005.

\bibitem{SVV} A. Simis, W. Vasconcelos, and R. Villarreal, On the ideal theory of graphs. \emph{J. Algebra} {\bf 167} (1994), 389--416.

\bibitem{Stanley} R. P. Stanley, Combinatorics and commutative algebra. Second edition. \emph{Progress in Mathematics} {\bf 41}. Birkh�ser Boston, Inc., Boston, MA, 1996.

\bibitem{V1} R. H. Villarreal,
Cohen-Macaulay graphs.
Manuscripta Math. {\bf 66} (1990), no. 3, 277--293.

\bibitem{Z} X. Zheng, Resolutions of facet ideals. Comm. Algebra {\bf 32}
(2004) 2301-2324.

\end{thebibliography}
\end{document}